\documentstyle[11pt]{article}
\title{A derivative formula for the free energy function
} 
\author{Yu Zhang}
\date{}
\begin{document}
\baselineskip .20in
\maketitle
\begin{abstract}
We consider bond  percolation on  the ${\bf Z}^d$ lattice. Let $M_n$ be the 
number of open clusters in $B(n)=[-n, n]^d$. It is well known that $E_pM_n / (2n+1)^d$  converges to 
the free energy function $\kappa(p)$ at the zero field.  
In this paper, we show that
$\sigma^2_p(M_n)/(2n+1)^d$ converges to $-(p^2(1-p)+p(1-p)^2)\kappa'(p)$. 
\end{abstract}

\section{ Introduction and statement of results.}

Consider bond percolation on  the ${\bf Z}^d$ lattice, 
in which  bonds are independently open with 
probability $p$ and closed with probability
$1-p$. 
The corresponding probability measure on the configurations of open
and closed bonds is denoted by $P_p$. We also denote by $E_p(X)$ 
and $\sigma_p^2(X)$ the expectation  and the variance of $X$ with respect to $P_p$. 
The open cluster of the vertex $x$, 
${\bf C}(x)$,
consists of all vertices that are connected to $x$ by an open path.
Here an open  path from $u$ to $v$ is a sequence $(v_0, b_0, v_1, ...,v_{i}, b_i, v_{i+1},...,v_n)$
with distinct vertices $v_i$ ($0\leq i\leq n$) and open bonds $b_i$ adjacent $v_i$ and $v_{i+1}$
such that $v_0=u$ and $v_n=v$. 
For vertex set $A$, $A_e$ denotes the bonds with  both  vertices in $A$.
Also, $|A|$  
denotes the cardinality of $A$, and $|A_e|$ denotes the number of bonds in $A_e$. 
We choose ${\bf 0}$ as the origin. The percolation probability is
\begin{eqnarray*}
\theta (p)= P_p(|{\bf C}({\bf 0})|=\infty),
\end{eqnarray*}
and the  critical probability is 
$$p_c=\sup\{p:\theta (p)=0\}.$$
We denote the open cluster distribution by
$$\theta_n(p)=P_p(|{\bf C}({\bf 0})|=n).$$
By analogy with the Ising model, we introduce the magnetization function as
$$M(p,h)=1-\sum_{n=0}^\infty \theta_n(p) e^{-nh}\mbox{ for }h\geq 0.$$
By setting $h=0$ in the magnetization function,
$$M(p,0)=\theta(p).$$
Using  term by term differentiation,
we also have
$$\lim_{h\rightarrow 0^+}{\partial M(p,h)\over \partial h}=E_p( |{\bf C}({\bf 0})|;|{\bf C}({\bf 0})|<\infty)=\chi^f(p).$$
$\chi^f(p)$ is called the {\em mean cluster size}. The {\em free energy} $F(p,h)$ is defined by
$$F(p,h)=h(1-\theta_0(p))+\sum_{n=1}^\infty {1\over n} \theta_n(p)e^{-hn}\mbox{ for } h >0.$$
If we differentiate with respect to $h$, then we find
$$ {\partial F(p,h)\over \partial h}=M(p,h).$$
For $h>0$, the free energy is infinitely differentiable with respect to $p$.
If $h=0$, $F(p, 0)$ is called the {\em zero-field} free energy.
The zero-field free energy $F(p,0)$ is a more interesting and more difficult object of study
since it is believed that there is a singularity point at $p_c$. By our definition,
$$F(p,0)=E(|{\bf C}({\bf 0})|^{-1}; |{\bf C}({\bf 0})| >0) .$$
Grimmett (1981) discovered that the zero-field  free energy  also coincides  with the number of open
clusters per vertex. 
Let us define the number of open clusters per vertex as follows.
We denote by  $M_n$ the 
number of open clusters in $B(n)$. By a standard ergodic theorem (see Theorem 4.2 in
Grimmett  (1989)), the limit
$$\lim_{n\rightarrow \infty}{1\over {|B(n)|}}M_n=\lim_{n\rightarrow \infty}{1\over {(2n+1)^d}}M_n=\kappa(p) 
\mbox { a.s. and }L_1\eqno{(1.1)}$$ 
exists for all $0\leq p\leq 1$.
$\kappa (p)$ is called the {\em number of  open clusters per vertex}. Grimmett (1981) proved  
that 
$$\kappa(p)=F(p,0). \eqno{(1.2)}$$
$\kappa(p)$, as a function of $p$,  
is   analytic for $p\neq p_c$ and
differentiable on $[0,1]$ (see Kesten (1982)). In particular, 
$\kappa(p)$ is proved  (see Kesten (1982)) to be twice differentiable at $p_c$ for the square lattice.  
In general,
physicists believe that the zero-field free energy is twice, but not three times  differentiable at $p_c$. 

On the other hand, for $\kappa(p)$, as a limit of random variables, 
Zhang (2001) showed the following central limit theorem. For $p\in [0, 1]$,
$${M_n - E_pM_n\over \sigma_p(M_n) }\Rightarrow \mbox{ a standard normal distribution.}\eqno{(1.3)}$$
Zhang also showed  large deviations for $M_n$. In this paper, we show another property for $\kappa(p)$.\\

{\bf Theorem.} {\em For $0\leq p\leq 1$,} 
$$\lim_{n\rightarrow \infty} \sigma^2_p(M_n)/(2n+1)^d =-(p^2(1-p)+p(1-p)^2)\kappa'(p).$$

{\bf Remark.} 
If $b_0$ is the bond with vertices ${\bf 0}$ and $(1,0,\cdots, 0)$, and
 ${\cal G}(b)$ is the event that there does not exist an open path in ${\bf Z}^d_e\setminus b$
 connecting the two vertices
of bond $b$, then by (2.9) in the proof of the theorem, we have
$$\kappa'(p)=-dP_p({\cal G}(b_0)).\eqno{(1.5)}$$
By using Sykes and Essam's formula (1964), we know that for the square lattice,
$$\kappa'(0.5)=-1.$$
Thus
$$\lim_{n\rightarrow \infty} \sigma^2_{0.5}(M_n)/(2n+1)^d =0.25.\eqno{(1.6)}$$

\section{ Proof of theorem.}

For any bond $b$, let $v_1(b)$ and $v_2(b)$
be the two vertices of $b$. 
Given $p$, we start by taking the  
derivative of $E_p(M_n)$. Note that
$$E_pM_n=\sum^\infty _{l=1}P_p(M_n\geq l),$$ and the
event $\{M_n\geq l\}$ is  
decreasing. Let $\{M_n\geq l\}(b)$ be the event that
$b$ is a pivotal bond for $\{M_n \geq l\}$. By Russo's formula, note that $P_p(M_n\geq 1)=1$, so
\begin{eqnarray*}
{d E_pM_n \over dp}
&=& -\sum_{l=2}^{\infty}\sum_{b\in B_e(n)} P_p(\{ M_n\geq l\}(b))\\
&=& -\sum_{b\in B_e(n)}\sum_{l=2}^\infty P_p(M_n=l-1\mbox
{ or $l$ if $b$ is open or closed})\\
&=&-\sum_{b\in B_e(n)}P_p\left(\bigcup_{l=2}^\infty M_n=l-1\mbox
{ or $l$ if $b$ is open or closed}\right). \hskip 4cm (2.1)
\end{eqnarray*}
Let
$${\cal E}_n(b)=\{b \mbox{ is a pivotal bond for the open connection
of $v_1(b)$ and $v_2(b)$ in $B(n)$}\}.$$
In other words, if $b$ is open, then $v_1(b)$ and $v_2(b)$ are connected
by open paths. Conversely, if  $b$ is closed, then $v_1(b)$ and $v_2(b)$ are not connected
by open paths. 
Thus, for each $b\in B_e(n)$,
$$P_p\left(\bigcup_{l=2}^\infty M_n=l-1\mbox
{ or $l$ if $b$ is open or closed}\right)
=P_p({\cal E}_n(b)).\eqno{(2.2)}$$ 

Let ${\cal G}_n(b)$ be the event that there does not exist an open path connecting $v_1(b)$ to $v_2(b)$ inside
$B_e(n)\setminus b$.  Then we would   have
$${\cal E}_n(b)={\cal G}_n(b).\eqno{(2.3)}$$
To see (2.3), if there were such a path, then $v_1(b)$ and $v_2(b)$ would always be
connected by an open path whenever $b$ is open or closed. 
So $b$ would not be a pivotal bond. 
On the other hand, if there
does not exist an open path connecting $v_1(b)$ and $v_2(b)$ in $B_e(n)\setminus b$,
then $b$ should be a pivotal bond for the open connection of $v_1(b)$ and $v_2(b)$.
\begin{figure}
\begin{center}
\setlength{\unitlength}{0.0125in}%
\begin{picture}(200,320)(67,750)
\thicklines
\put(100,770){\framebox(200,200)[br]{\mbox{$[-n,n]^2$}}}
\put(175,880){\mbox{$b$}{}}
\put(170, 890){\line(-1,0){10}}
\put(170, 891){\line(-1,0){10}}
\put(170,890){\circle*{4}}
\put(160,890){\circle*{4}}
\put(166,885){\circle*{3}}
\put(166,880){\circle*{3}}
\put(166,875){\circle*{3}}
\put(166,870){\circle*{3}}
\put(166,865){\circle*{3}}
\put(166,860){\circle*{3}}
\put(171,860){\circle*{3}}
\put(175,860){\circle*{3}}
\put(175,855){\circle*{3}}
\put(175,850){\circle*{3}}
\put(175,845){\circle*{3}}
\put(175,840){\circle*{3}}
\put(175,835){\circle*{3}}
\put(175,830){\circle*{3}}
\put(175,825){\circle*{3}}
\put(175,820){\circle*{3}}
\put(175,815){\circle*{3}}
\put(175,810){\circle*{3}}
\put(175,805){\circle*{3}}
\put(175,800){\circle*{3}}
\put(175,795){\circle*{3}}
\put(175,790){\circle*{3}}
\put(175,785){\circle*{3}}
\put(175,780){\circle*{3}}
\put(175,775){\circle*{3}}
\put(175,770){\circle*{3}}

\put(166,895){\circle*{3}}
\put(166,900){\circle*{3}}
\put(166,905){\circle*{3}}
\put(160,905){\circle*{3}}
\put(155,905){\circle*{3}}
\put(150,905){\circle*{3}}
\put(145,905){\circle*{3}}
\put(140,905){\circle*{3}}
\put(135,905){\circle*{3}}
\put(130,905){\circle*{3}}
\put(125,905){\circle*{3}}
\put(120,905){\circle*{3}}
\put(115,905){\circle*{3}}
\put(110,905){\circle*{3}}
\put(105,905){\circle*{3}}
\put(100,905){\circle*{3}}

\put(160, 890){\line(0,-1){50}}
\put(160, 840){\line(-1,0){80}}
\put(80, 840){\line(0,1){90}}
\put(80, 930){\line(-1,0){40}}
\put(40, 930){\line(0,1){80}}
\put(40, 1010){\line(1,0){80}}
\put(120, 1010){\line(0,1){20}}
\put(120, 1030){\line(1,0){50}}
\put(170, 1030){\line(0,-1){140}}

\end{picture}
\end{center}
\caption{ The figure shows that ${\cal G}_n(b)$ occurs, but ${\cal G}(b)$ does not occur. The two
vertices of $b$ are connected by  open paths, but they have to reach to the boundary of $B(n)$
before connecting with each other. The dotted path indicates the closed bonds that block the connection of open paths
from $v_1(b)$ and $v_2(b)$ inside $B(n)$.}
\end{figure}
With these observations, 
$$
{d E_pM_n \over dp}= - \sum_{b\in B_e(n)} P_p({\cal G}_n(b))=- \sum_{b\in B(n)}P_p({\cal G}(b))-\sum_{b\in B_e(n)} P_p\left({\cal G}_n(b)\cap {\cal G}^C(b)\right).
$$
By translation invariance, the ratio of the first term above is
$$\sum_{b\in B_e(n)}P_p({\cal G}(b))/|B_e(n)|=P_p({\cal G}(b_0)).\eqno{(2.4)}$$
We use the following lemma to estimate the second term. \\

{\bf Lemma.} 
$$\lim_{n\rightarrow \infty}{1\over |B_e(n)|}
\sum_{b\in B_e(n)} P_p\left({\cal G}_n(b)\cap {\cal G}^C(b)\right)=0\mbox{
{\em uniformly on }$[0,1]$.}$$

{\bf Proof.} Let $A(n, m)= B(n)\setminus B(m)$ for $m\leq n$, and let $\bar{A}(n,m)$ be
the closure of $A(n,m)$. Then
\begin{eqnarray*}
&&{1\over |B_e(n)|}\sum_{b\in B_e(n)} P_p\left({\cal G}_n(b)\cap {\cal G}^C(b)\right)\\
&\leq &{1\over |B_e(n)|}\sum_{b\in B_e(n-\sqrt{n})} P_p\left({\cal G}_n(b)\cap  {\cal G}^C(b)\right)
+ {1\over |B_e(n)|}\sum_{b\in \bar{A}_e(n, n-\sqrt{n})} P_p\left({\cal G}_n(b)\cap {\cal G}^C(b)\right)\\
&\leq & {1\over |B_e(n)|}\sum_{b\in B_e(n-\sqrt{n})} P_p\left({\cal G}_n(b)\cap {\cal G}^C(b)\right)+O\left(n^{-0.5}\right),\hskip 6cm (2.5)
\end{eqnarray*}
where we may assume that $n-\sqrt{n}$ is an integer; otherwise we may use $\lceil n-\sqrt{n} \,\,\rceil$ to
replace $n-\sqrt{n}$.
Let us estimate the first term in the above inequality.
For $b\in B_e( n-\sqrt{n})$,  ${\cal E}_n(b)\cap {\cal E}^C(b)$ implies that
there exists an open path from $v_1(b)$ and $v_2(b)$ without using $b$, but open paths
cannot stay inside $B(n)$. In other words, 
the open
path adjacent to $v_1(b)$ has to reach  the boundary of $B(n)$ before reaching  $v_2(b)$ (see Fig. 1). 
Similarly, the open
path adjacent to $v_2(b)$ has to reach  the boundary of $B(n)$ before reaching  $v_1(b)$ (see Fig. 1).
Therefore, there exist two disjoint open paths from $v_1(b)$ and $v_2(b)$ such that both reach
to the boundary of $v_1(b)+ B(\sqrt{n}-1)$.  Let ${\cal D}(b, \sqrt{n})$ be the event.
By using the estimate of Theorem 6.1 in Grimmett (1989), for all $p\in [0,1]$, there exist
constant $C=C(d)$ and  $\delta=\delta(d)\leq 0.5$
such that
$$P_p\left({\cal D}(b, \sqrt{n})\right)\leq Cn^{-\delta}.\eqno{(2.6)}$$
By (2.5) and (2.6),
$${1\over |B_e(n)|}\sum_{b\in B_e(n)} P_p\left({\cal G}(b)\cap{\cal G}^C_n(b)\right)\leq O(n^{-\delta}).\eqno{(2.7)}$$
Therefore, the lemma follows from (2.7). $\Box$\\

It follows from (2.3), the lemma, and (2.4) that
$$
{d E_pM_n \over |B(n)|dp} =- (|B_e(n)|/| B(n)|)\sum_{b\in B(n)}P_p({\cal G}(b))/|B_e(n)|+O(n^{-\delta}).\eqno{(2.8)}
$$
Note that 
$$|B_e(n)|=2dn(2n+1)^{d-1}\mbox{ and } |B(n)|=(2n+1)^d,$$
so if we let $n\rightarrow \infty$ in (2.8), 
$$\kappa'(p) = -\lim_{n\rightarrow \infty} \sum_{b\in B_e(n)}P_p({\cal G}_n(b))/(2n+1)^d=-dP_p\left({\cal G}(b_0)\right).\eqno{(2.9)}$$
Now we estimate the variance of $M_n$. We list the bonds in $B_e(n)$ in some order:
$$\{b_1, \cdots, b_k\}.$$
We define the independent Bernoulli-random variables  $\{\omega(b_i)\}$ 
for $1\leq i\leq k$ to be $\omega(b_i)=0$ (open) or $\omega(b_i)=1$ (closed) with probability $p$ or $1-p$.
Now we construct the following filtration:
$${\cal F}_0=\{\emptyset, \Omega\}\subset {\cal F}_1=\{\sigma-\mbox{field generated by }\omega(e_1)\}\subset,\cdots,
\subset {\cal F}_k=\{\sigma-\mbox{field generated by }\omega(e_1),\cdots, \omega(e_k)\}.$$
The martingale representation of $M_n-E_pM_n$ is
$$M_n-E_p(M_n)=\sum_{i=0}^{k-1} [E_p(M_n|{\cal F}_{i+1})-E_p(M_n|{\cal F}_{i})].$$
Let the martingale difference be
$$\Delta_{i, k}=[E_p(M_n|{\cal F}_{i+1})-E_p(M_n|{\cal F}_{i})].$$
The variance is
$$\sigma^2_p(M_n)=\sum_{i=0}^{k-1} E_p (\Delta_{i, k}^2).\eqno{(2.10)}$$
Both $M_n$ and $\Delta_{i,k}$ can be viewed as functions on $[0, 1]^k$ and $[0,1]^{i}$, respectively. So we can write
$M_n(c_1, \cdots, c_k)$ and $\Delta_{i, k}(c_1,\cdots, c_i)$ for them, where
$c_i$ only takes a value of zero or one. Thus
\begin{eqnarray*}
&&\Delta_{i,k}=\Delta_{i,k}(c_{1},\cdots, c_i)\\
&=&\sum_{c_{i+1},\cdots, c_k} M_n(c_{1},\cdots, c_i, c_{i+1},\cdots, c_k)P_p\left(\omega(b_{i+1})=c_{i+1}, \cdots, \omega(b_{k})=c_{k}\right)\\
&&-\sum_{c'_i, c_{i+1},\cdots, c_k} M_n(c_{1},\cdots, c_{i-1},c_i', c_{i+1},\cdots, c_k)P_p(\omega(b_{i})=c_{i}',\omega(b_{i+1})=c_{i+1}, \cdots, \omega(b_{k})=c_{k})\\
&=&\sum_{c'_i, c_{i+1},\cdots, c_k}[M_n(c_{1},\cdots, c_i, \cdots, c_k)-M_n(c_{1},\cdots, c_{i-1},c_i',c_{i+1}, \cdots, c_k)]\\
&&\,\,\,\,\,\,\,\,\cdot P_p(\omega(b_{i})=c_{i}', \omega(b_{i+1})=c_{i+1},\cdots, \omega(b_{k})=c_{k}),\hskip 6cm (2.11)
\end{eqnarray*}
where the sum takes over all possible values
of $c_i$ and $c_i'$. 
On ${\cal G}_n^C(b_i)$,
\begin{eqnarray*}
&&M_n(c_{i+1},\cdots, c_i, \cdots, c_k)-M_n(c_{i+1},\cdots, c_i'=0, \cdots, c_k)=0 \mbox{ for $c_i=0$ or $c_i=1$},\\
&&M_n(c_{i+1},\cdots, c_i, \cdots, c_k)-M_n(c_{i+1},\cdots, c_i'=1, \cdots, c_k)=0 \mbox{ for $c_i=0$ or $c_i=1$}.
\end{eqnarray*}
Thus, by (2.11), 
$$\Delta_{i,k}\left(1-I_{{\cal G}_n(b_i)}\right)=0,\eqno{(2.12)}$$
where $I_{\cal A}$ is the indicator of ${\cal A}$.
 With this observation,
$$\Delta_{i,k}=\Delta_{i,k}I_{{\cal G}_n(b_i)}.\eqno{(2.13)}$$
Note that on ${{\cal G}_n(b_i)}$,
$$[M_n(c_{1},\cdots, c_i=1, \cdots, c_k)-M_n(c_{1},\cdots, c_{i-1}, c_i'=0, c_{i+1},\cdots, c_k)]=-1.$$
so
\begin{eqnarray*}
&&\Delta_{i,k}(c_1, \cdots, c_{i-1}, 1)I_{{\cal G}_n(b_i)}\\
&=&I_{{\cal G}_n(b_i)}\sum_{ c_{i+1},\cdots, c_k}[M_n(c_{1},\cdots, c_i=1, \cdots, c_k)-M_n(c_{1},\cdots, c_{i-1}, c_i'=0, c_{i+1},\cdots, c_k)]\\
&& \,\,\,\,\,\,\,\,\,\cdot P_p(\omega(b_i)=c_i'=0, \omega(b_{i+1})= c_{i+1},\cdots, \omega(b_{k})=c_{k})\\
&&+I_{{\cal G}_n(b_i)}\sum_{ c_{i+1},\cdots, c_k}[M_n(c_{1},\cdots, c_i=1, \cdots, c_k)-M_n(c_{1},\cdots, c_{i-1}, c_i'=1, c_{i+1},\cdots, c_k)]\\
&& \,\,\,\,\,\,\,\,\,\cdot P_p(\omega(b_i)=c_i'=1, \omega(b_{i+1})= c_{i+1},\cdots, \omega(b_{k})=c_{k})\\
&=&I_{{\cal G}_n(b_i)}\sum_{ c_{i+1},\cdots, c_k}[M_n(c_{1},\cdots, c_i=1, \cdots, c_k)-M_n(c_{1},\cdots, c_{i-1}, c_i'=0, c_{i+1},\cdots, c_k)]\\
&& \,\,\,\,\,\,\,\,\,\cdot P_p(\omega(b_i)=c_i'=0)P_p( \omega(b_{i+1})= c_{i+1},\cdots, \omega(b_{k})=c_{k})\\
&=&-(1-p)I_{{\cal G}_n(b_i)}
\end{eqnarray*}
Therefore,
$$\Delta_{i,k}(c_1, \cdots, c_{i-1}, 1)I_{{\cal G}_n(b_i)}=(1-p)^2I_{{\cal G}_n(b_i)}.\eqno{(2.14)}$$

Similarly, 
note that on ${{\cal G}_n(b_i)}$,
$$[M_n(c_{1},\cdots, c_i=0, \cdots, c_k)-M_n(c_{1},\cdots, c_{i-1}, c_i'=1, c_{i+1},\cdots, c_k)]=1.$$
so
\begin{eqnarray*}
&&\Delta_{i,k}(c_1, \cdots, c_{i-1}, 0)I_{{\cal G}_n(b_i)}\\
&=&I_{{\cal G}_n(b_i)}\sum_{ c_{i+1},\cdots, c_k}[M_n(c_{1},\cdots, c_i=0, \cdots, c_k)-M_n(c_{1},\cdots, c_{i-1}, c_i'=0, c_{i+1},\cdots, c_k)]\\
&& \,\,\,\,\,\,\,\,\,\cdot P_p(\omega(b_i)=c_i'=0, \omega(b_{i+1})= c_{i+1},\cdots, \omega(b_{k})=c_{k})\\
&&+I_{{\cal G}_n(b_i)}\sum_{ c_{i+1},\cdots, c_k}[M_n(c_{1},\cdots, c_i=0, \cdots, c_k)-M_n(c_{1},\cdots, c_{i-1}, c_i'=1, c_{i+1},\cdots, c_k)]\\
&& \,\,\,\,\,\,\,\,\,\cdot P_p(\omega(b_i)=c_i'=1, \omega(b_{i+1})= c_{i+1},\cdots, \omega(b_{k})=c_{k})\\
&=&I_{{\cal G}_n(b_i)}\sum_{ c_{i+1},\cdots, c_k}[M_n(c_{1},\cdots, c_i=0, \cdots, c_k)-M_n(c_{1},\cdots, c_{i-1}, c_i'=1, c_{i+1},\cdots, c_k)]\\
&& \,\,\,\,\,\,\,\,\,\cdot P_p(\omega(b_i)=c_i'=1)P_p( \omega(b_{i+1})= c_{i+1},\cdots, \omega(b_{k})=c_{k})\\
&=&pI_{{\cal G}_n(b_i)}
\end{eqnarray*}
Therefore,
$$\Delta_{i,k}(c_1, \cdots, c_{i-1}, 0)I_{{\cal G}_n(b_i)}=p^2I_{{\cal G}_n(b_i)}.\eqno{(2.15)}$$

Together with (2.13), (2.14) and (2.15),
\begin{eqnarray*}
E_p\Delta_{i,k}^2&=&E_p\Delta_{i,k}^2I_{{\cal G}_n(b_i)}\\
&=&E_p \Delta_{i,k}^2(c_1, \cdots, c_{i-1}, 1)I_{\{\omega(b_i)=1\}}I_{{\cal G}_n(b_i)}+
E_p \Delta_{i,k}^2(c_1, \cdots, c_{i-1}, 0)I_{\{\omega(b_i)=0\}}I_{{\cal G}_n(b_i)}\\
&=&E_p(1-p)^2I_{\{\omega(b_i)=1\}}I_{{\cal G}_n(b_i)}+E_pp^2I_{\{\omega(b_i)=0\}}I_{{\cal G}_n(b_i)}.
\end{eqnarray*}
Note that $\{\omega(b_i)\}$ and $I_{{\cal G}_n(b_i)}$ are independent, so 
$$E_p\Delta_{i,k}^2=[(1-p)^2 p+ p^2(1-p)]E_pI_{{\cal G}_n(b_i)}=[(1-p)^2 p+ p^2(1-p)]P_p({\cal G}_n(b_i)).\eqno{(2.16)}$$
By  (2.10) and (2.16),
$$
\sigma^2_p(M_n)=\sum_{i=1}^k  E_p (\Delta_{i, k}^2)= [(1-p)^2 p+ p^2(1-p)]\sum_{b\in B_e(n)}P_p({\cal G}_n(b)). \eqno{(2.17)}$$
Therefore,  if we divide the both sides of (2.17) by $(2n+1)^d$ and let $n$ go to infinity, by (2.9),
\begin{eqnarray*}
\lim_{n\rightarrow \infty}\sigma^2_p(M_n)/(2n+1)^d&=&\lim_{n\rightarrow \infty}[(1-p)^2 p+ p^2(1-p)]\sum_{b\in B_e(n)}P_p({\cal G}_n(b))/(2n+1)^d\\
&=&-[(1-p)^2 p+ p^2(1-p)]\kappa'(p).
\end{eqnarray*}
So the theorem follows.

\begin{center}{\large \bf References} \end{center}
Grimmett, G. (1981). On the differentiability of the number of clusters per vertex in the percolation model. 
{\em J. Lond. Math. Soc.} {\bf 23} 372--384.\\
Grimmett, G. (1989). {\em Percolation}. Springer-Verlag, New York.\\
Kesten, H. (1982). {\em Percolation Theory for Mathematicians}. Birkhauser, Boston.\\
Sykes, M. F., and Essam, J. W. (1964). Exact critical percolation probabilities for site and bond problems in two dimensions. {\em J. Math. Phys.} {\bf 5} 1117-1127.\\
Zhang, Y. (2001). A martingale approach in the study of percolation clusters on the $Z^d$ lattice. 
{\em J. Theore. Probab.} {\bf 14} 165--187.\\

Yu Zhang\\
Department of Mathematics\\
University of Colorado\\
Colorado Springs, CO 80933\\
yzhang3@uccs.edu

\end{document}